\documentclass[11pt,twoside,reqno]{amsart}
\usepackage[all]{xy}
\usepackage{amssymb,latexsym,amsthm,amsmath,mathtools}
\usepackage{enumitem,color}

\usepackage[colorlinks=true, linkcolor=black, citecolor=blue]{hyperref}
\usepackage[capitalize,nameinlink]{cleveref}
\usepackage{cleveref}
\usepackage[
    backend=biber,
    style=numeric-comp,
    doi=true,
    url=false,
    maxbibnames=99,
    giveninits=true
]{biblatex}
\DeclareFieldFormat[article]{title}{#1}
\DeclareFieldFormat{title}{#1}
\DeclareFieldFormat[incollection]{title}{#1}
\crefname{section}{Section}{Sections}
\crefname{subsection}{Subsection}{Subsections}
\crefname{subsubsection}{Subsubsection}{Subsubsections}
\usepackage{parskip}
\usepackage{mathtools}

\usepackage[backend=biber, doi=true, url=false, maxbibnames=99, giveninits=true]{biblatex}
\usepackage{orcidlink}
\long\def\symbolfootnote[#1]#2{\begingroup%
\usepackage{hyperref}
\usepackage{tikz}
\usepackage{blkarray}
\usepackage{arydshln}
\usepackage{tabularx}
\setlength\dashlinedash{0.5pt}
\setlength\dashlinegap{1.5pt}
\setlength\arrayrulewidth{0.3pt}
\def\thefootnote{\fnsymbol{footnote}}\footnote[#1]{#2}\endgroup}

\makeatletter
\def\imod#1{\allowbreak\mkern10mu({\operator@font mod}\,\,#1)}
\makeatother

\newtheorem{theorem}{Theorem}[section]
\newtheorem{lemma}[theorem]{Lemma}
\newtheorem{corollary}[theorem]{Corollary}
\newtheorem{proposition}[theorem]{Proposition}
\newtheorem*{theorem*}{Theorem}
\theoremstyle{definition}

\newtheorem{question}[theorem]{Question}

\numberwithin{equation}{section}

\newcommand{\ignore}[1]{}

\newcommand{\mynote}[1]{}
\begin{document} 
\author[Saini P.]{Prachi Saini\,\orcidlink{0009-0006-7627-1538}}
\address{Indian Institute of Science Education and Research Mohali, Knowledge City, Sector 81, Punjab 140306, India}
\email{prachi2608saini@gmail.com}

\thanks{The author acknowledges financial support from the Indian Institute of Science Education and Research (IISER) Mohali through a Postdoctoral Fellowship.}

\setcounter{section}{0}
\title{Sum of Powers with Matrix Coefficients}


\begin{abstract}
Let $F$ be an algebraically closed field of characteristic
zero. We consider polynomial maps with constants
$\omega(x,y)=A_1x^k+A_2y^k$ on $M_n(F)$, where $A_1$
is invertible and $k\geq2$. We show that the necessary condition for surjectivity
obtained by Saini and Singh
(Polynomial maps with constants on matrix algebra,
Linear Algebra Appl. 752 (2027), 256--277)
is not sufficient in any dimension $n\geq6$. We also completely characterize surjectivity for $n=5$ and $6$. 
\end{abstract}

\maketitle


\section{Introduction}

Let $F$ be a field and let $A$ be an associative $F$-algebra. A polynomial $f\in F\langle x,y,\ldots, x_m\rangle$ in non-commuting variables defines a map on $A$ by evaluation. If $f$ is a polynomial identity for $A$, then its image is $\{0\}$. Otherwise, one may ask which elements of $A$ belong to its image. If the map is not surjective, can every element of $A$ be expressed as a sum of a prescribed number of values of $f$?

A central direction concerns Waring-type problems:
whether every element of an algebra can be expressed
as a sum of a prescribed number of polynomial values.
Representations of matrices as sums of powers form
an important class of such problems; see
\cite{BresarSemrlwaringproblem, MR4575412, MR4592319, MR4526834}.
Related questions concern the structure of polynomial
images themselves, including the L'vov--Kaplansky
conjecture on images of multilinear polynomials.
For an account of polynomial images and
Waring-type problems, we refer to
\cite{bresar2026waring}.

A further variation is obtained by allowing the
coefficients of a polynomial to belong to the algebra
itself. These polynomials are called generalized
polynomials, or polynomials with constants. In this
setting, the image depends on the chosen coefficients
as well as on the polynomial expression; see \cite{demello2026generalizedcommutatoroctonionalgebra, PolynomialswithconstantsM, PolynomialsonSplitOctonion, SAINI2027256}. Moreover,
the image need not be invariant under conjugation,
since conjugating a value also conjugates the
coefficients. The study of such maps arises naturally as analogues of word maps with constants in group theory; see \cite{GnutovWordmapswithcosntants, WordmapsandWaringproblemLArsen, Wordmapswithconstants}

Panja, Saini and Singh \cite{PolynomialswithconstantsM} studied
polynomial maps with constants on $M_2(F)$ over an
algebraically closed field. In particular, they
characterized the surjectivity of maps of the form
$A_1x^{k_1}+A_2y^{k_2}$ and proved that the image of
$A_1xy-A_2yx$ is a vector space. Subsequently, Saini
and Singh \cite{SAINI2027256} considered the maps
$$
\omega\colon M_n(F)\times M_n(F)\longrightarrow M_n(F),
\qquad
\omega(x,y)=A_1x^k+A_2y^k,
$$
where $F$ is algebraically closed of characteristic
zero, $k\geq2$, and $A_1$ is invertible. Writing $r_0=\operatorname{nullity}(A_2)$, they proved
that $n>k(r_0-1)$
is necessary for surjectivity and that this condition
is also sufficient for $n=3,4$.
This raises the question of whether the same
inequality characterizes surjectivity in higher
orders.

We prove that the necessary inequality is sufficient
for $n=5$, but not for $n=6$.
More generally, for every $n\geq6$ and $k\geq n-1$,
we characterize surjectivity in terms of the nullity
of the second coefficient.
Our main results are the following.
Throughout, $F$ is algebraically closed of
characteristic zero.
\begin{theorem}
\label{thm:large-exponents}
Let $F$ be an algebraically closed field of characteristic zero. Let $n\geq6$ and $k\geq n-1$, and
$A_1,A_2\in M_n(F)$ with $A_1$ invertible.
Then the map
$$
\omega\colon M_n(F)\times M_n(F)\longrightarrow M_n(F),
\qquad
\omega(x,y)=A_1x^k+A_2y^k,
$$
is surjective if and only if
$\operatorname{nullity}(A_2)\leq1$.
\end{theorem}

\begin{theorem}
\label{thm:n5-n6}
Let $F$ be an algebraically closed field of characteristic
zero. Let $n\in\{5,6\}$, and $A_1,A_2\in M_n(F)$
with $A_1$ invertible. For an integer $k\geq2$, consider
the map
\[
\omega\colon M_n(F)\times M_n(F)\longrightarrow M_n(F),
\qquad
\omega(x,y)=A_1x^k+A_2y^k.
\]
Set $r_0=\operatorname{nullity}(A_2)$. Then:
\begin{enumerate}
    \item  For $n=5$, the map $\omega$ is surjective
    if and only if $5>k(r_0-1)$.

    \item For $n=6$, the map $\omega$ is surjective
    if and only if $6>k(r_0-1)$ and
    $(r_0,k)\neq(2,5)$.
\end{enumerate}
\end{theorem}

The paper is organized as follows. Section~\ref{sec: Prelims} contains the reductions and auxiliary lemmas used throughout. In Section~\ref{sec: not surjective}, we establish the general obstruction and deduce the classification for large exponents. Sections~\ref{sec:n5} and~\ref{sec:n6} give the classifications in orders five and six. In Section~\ref{sec: concluding}, we conclude with a question concerning higher orders.

\section{Preliminaries}
\label{sec: Prelims}
Throughout, we assume $F$ is an algebraically closed field of characteristic $0$ and $k\geq 2$. We recall the following reductions \cite{SAINI2027256}. Let $A_1, A_2 \in M_n(\mathbb{F})$ be non-zero elements such that one of the coefficients is invertible. Consider the map
$$\omega : M_n(\mathbb{F})^2 \to M_n(\mathbb{F}), \qquad
\omega(x,y) = A_1x^k + A_2y^k.$$

Without loss of generality, we assume $A_1$ is invertible. Define $B = A_1^{-1}A_2 \in M_n(\mathbb{F})$ and $r_0=\operatorname{nullity}(B)=\operatorname{nullity}(A_2)$.  

Moreover, the map $\omega$ is surjective if and only if the map $\omega_B(x,y) = x^{k} + B y^{k}$ is surjective. Indeed, $\omega(x,y)=C$ if and only if $\omega_B(x,y)=A_1^{-1}C$.

\textbf{Reduction to Jordan form.}
Surjectivity of the map $\omega$ remains invariant under the action of the automorphism group of $\mathrm{PGL}_n(F)$. Thus, we may replace $B$ by its Jordan canonical form
$$
J = \left(\bigoplus_{i=1}^{r'} J_{\lambda_i,n_i}\right)\oplus \left(\bigoplus_{s=1}^{r_0} J_{0,m_s}\right),
$$
where $\lambda_i\neq 0$ and $\sum_{s=1}^{r_0} m_s = n_0$. For $1\le s\le r_0$, set $\ell_s := n-n_0+m_1+\cdots+m_s$, and define $\sigma:\{1,\dots,n\}\to\{1,\dots,n\}$ by
$$
\sigma(i)=
\begin{cases}
i, & 1\le i\le n-n_0,\\
i-\#\{\,s:\ell_s<i\,\}, & n-n_0<i\le n,\ i\notin\{\ell_1,\dots,\ell_{r_0}\},\\
n-r_0+s, & i=\ell_s
\end{cases}
$$
associated to $J$ and $P$ be the permutation matrix corresponding to $\sigma$. Consider the partition of $C$ and $J$ obtained under the conjugation of $P$:
$$
 P^{-1}CP = \begin{pmatrix}
C_{11}^P & C_{12}^P\\
C_{21}^P & C_{22}^P
\end{pmatrix}\qquad \text{and} \qquad P^{-1}JP=\begin{pmatrix}
    \widehat{J}\\ 0
\end{pmatrix},
$$ 
where $C_{11}^P\in M_{n-r_0}(F)$, $C_{22}^P\in M_{r_0}(F)$ and $\widehat{J}\in M_{(n-r_0)\times n}(F)$. Note that $\widehat{J}$ is a full row rank. We set $L_1=\begin{pmatrix}
    C_{11}^P & C_{12}^P
\end{pmatrix}\in M_{n-r_0\times n}(F)$ and $L_2=\begin{pmatrix}
    C_{21}^P & C_{22}^P
\end{pmatrix}\in M_{r_0\times n}(F)$.

Denote $W_1\subseteq F^n$ as the subspace spanned by the first $n-r_0$ rows of $C^P$ and $W_2\subseteq F^n$ as the subspace spanned by the last $r_0$ rows of $C^P$.
In the following sections, we consider $C$ and $B$ up to conjugation by $P$.
\subsection{Notation}
Throughout, $F$ denotes an algebraically closed field of characteristic $0$, and the split central simple algebra over $F$ is $M_n(F)$, the algebra of all $n \times n$ matrices. We denote the automorphism group of $M_n(F)$ by $\mathrm{PGL}_n(F)$, the projective general linear group, which acts on $M_n(F)$ by conjugation.

For a matrix $C \in M_n(F)$, we denote its rows by $c_1, \dots, c_n \in F^n$ and $C=(c_1, \ldots, c_n)^\top$. For $1 \leq i \leq n$, let $e_i^\top$ denote the $i$-th standard basis vector of $F^n$. For $C_1 \in M_{n_1}(F)$ and $C_2 \in M_{n_2}(F)$, we define their direct sum by $ C_1 \oplus C_2 = \begin{pmatrix}
C_1 & 0 \\ 0 & C_2 \end{pmatrix} \in M_{n_1+n_2}(F)$. For $B \in M_n(F)$, we write $B \sim \bigoplus_{i=1}^r J_{\lambda_i,n_i}, \sum_{i=1}^r n_i = n$, to denote a Jordan canonical form of $B$, where $J_{\lambda_i,n_i}$ is the Jordan block of size $n_i$ corresponding to the eigenvalue $\lambda_i$. We denote by $r_0(B)$ the nullity of $B$, which turns out to be the number of Jordan blocks corresponding to the eigenvalue $0$, and by $r'(B)$ the number of Jordan blocks corresponding to nonzero eigenvalues. Thus, $ r(B) = r_0(B) + r'(B)$. 

For $S\in \mathrm{PGL}_n(F)$ and $C\in M_n(F)$, we denote $SCS^{-1}$ by $C^S$. For $Z\in M_{s\times t}(F)$, $\operatorname{span}(Z)$ is the span of rows of $Z$.
We retain the above notation throughout the paper. 

\subsection{Auxiliary Lemmas}
We first recall the Jordan form of the power of a nilpotent Jordan block. 
\begin{lemma}[Miller]
\label{lem:Miller}
Let $n > k \ge 2$. Then the $k$-th power of the nilpotent Jordan block $J_{0,n}$ is conjugate to
$$
\left(\bigoplus_{k-m} J_{0,\lfloor n/k \rfloor}\right)
\oplus
\left(\bigoplus_{m} J_{0,\lceil n/k \rceil}\right),
$$
where $m \equiv n \pmod{k}$ and $0 \le m \le k-1$.
\end{lemma}

\begin{lemma}\label{lem: rank-one}
For $n\geq 3$, every rank-one matrix in $M_n(F)$ is a $k$-th power
for $2\leq k\leq n-1$.
\end{lemma}
\begin{proof}
Let $Z$ have rank one. If $\operatorname{tr}(Z)\neq0$,
then $Z$ is similar to
$\operatorname{diag}(\operatorname{tr}(Z),0_{n-1})$
and hence has a $k$-th root. Otherwise, $Z$ is similar
to $J_{0,2}\oplus0_{n-2}$. By Miller's \cref{lem:Miller}, 
$$
J_{0,2}\oplus 0_{n-2}\sim \left(J_{0,k+1}\oplus0_{n-(k+1)}\right)^k.
$$
The zero summand is omitted when $k=n-1$. Therefore, $Z=X^k$ for some $X\in M_n(F)$.
\end{proof}
\begin{lemma}\label{lem: square completion}
For $1\leq s\leq n$, let $R\in M_{s\times n}(F)$. If $\operatorname{rank}(R)\leq n-s$ or $\operatorname{rank}(R)=s$,
then there exist $H\in M_{n-s\times n}(F)$ such that the matrix 
$\begin{pmatrix}
    H\\ R
\end{pmatrix}=X^2$ for some $X\in M_n(F)$.
\end{lemma}
\begin{proof}
If $\operatorname{rank}(R)=s$, extend its rows to an invertible matrix of order $n$ and take a square root.

Otherwise, set $t=\operatorname{rank}R\leq n-s$. For $t=0$, take $X=0$. Therefore, we take $t\neq 0$. Choose $Q\in\mathrm{GL}_s(F)$ such that the last $s-t$ rows of $Q^{-1}R$ are zero because $\operatorname{rank}(R)=t<s$. Let $S=I_{n-s}\oplus Q$ and write
$$
Q^{-1}RS=\{
u_1, \dots, u_t, 0_{s-t}\}^\top.
$$
Since $t\leq n-s$, define
$$
Z=
\{
u_1, \dots, u_t,
0_{(n-s-t)\times n},
e_1^\top, \dots, e_t^\top,
0_{(s-t)\times n}
\}^\top.
$$
The last $s$ rows of $ Z^2$ are $Q^{-1}RS$.
Let $X=SZS^{-1}$. Since $S=I_{n-s}\oplus Q$, the last
$s$ rows of $X^2=SZ^2S^{-1}$ are
$Q(Q^{-1}RS)S^{-1}=R$. Therefore, $ X^2$ has the required form.
\end{proof}
As a corollary, we obtain:
\begin{corollary}
    For $s\geq 1$, let $R\in M_{s\times n}(F)$. If $\operatorname{rank}(R)\leq n-s$ or $\operatorname{rank}(R)=s$,
then there exist $H\in M_{n-s\times n}(F)$ such that the matrix $\begin{pmatrix}
    R\\ H
\end{pmatrix}=X^2$ for some $X\in M_n(F)$.
\end{corollary}
\begin{proof}
Let $P=
\begin{pmatrix}
0_{s\times n-s}&I_s\\
I_{n-s}&0_{n-s\times s}
\end{pmatrix}$. Since $\operatorname{rank}(RP)=\operatorname{rank}R$, by
Lemma~\ref{lem: square completion} there exists $Z\in M_n(F)$
such that the last $s$ rows of $Z^2$ are $RP$. Then $(PZP^{-1})^2=PZ^2P^{-1}$ and hence the first $s$ rows of $(PZP^{-1})^2$ are therefore $RPP^{-1}=R$.
\end{proof}
\begin{lemma}\label{lem:residual lifting}
For $1\leq s\leq n$ and $2\leq k\leq n-1$, let $D\in M_{s\times n}(F)$ such that $\operatorname{rank}(D)$ is $s$.
If $H\in M_{s\times n}(F)$ satisfying $\operatorname{rank}(H)\in\{0,1,s\},$
then there exists $Y\in M_n(F)$ such that $$DY^k=H.$$
\end{lemma}
\begin{proof}
If $H=0$, take $Y=0$. Suppose that $\operatorname{rank}(H)=1$.
Write $H=bw$, where $b\in F^s$ is a nonzero column
vector and $w\in M_{1\times n}(F)$ is nonzero.
Since $D$ has full row rank, there exists $q\in F^n$
such that $Dq=b$. In particular, $q\neq0$, so $qw$
has rank one. By Lemma~\ref{lem: rank-one},
there exists $Y\in M_n(F)$ such that $Y^k=qw$.
Consequently, $$DY^k=D(qw)=(Dq)w=bw=H.$$

Finally, suppose that $\operatorname{rank}(H)=s$.
Extend the rows of $D$ and $H$ to bases of $F^n$,
obtaining invertible matrices
$$
S=\begin{pmatrix}D\\D_0\end{pmatrix},
\qquad
T=\begin{pmatrix}H\\H_0\end{pmatrix}.
$$
With $E=(I_s\ \ 0)$, we have $D=ES$ and $H=ET$.
The invertible matrix $S^{-1}T$ has a $k$-th root
$Y\in M_n(F)$. Hence
$$
DY^k=ES(S^{-1}T)=ET=H.
$$
When $s=n$, the additional rows are omitted.
\end{proof}
\begin{lemma}\label{lem:invertible-completion}
Let $1\leq r<n$, and set $s=n-r$.
Let $R\in M_{r\times n}(F)$ have rank $r$.
For every $L\in M_{s\times n}(F)$, there exists
$L'\in M_{s\times n}(F)$ such that
$
N=\begin{pmatrix}L'\\R\end{pmatrix}
$
is invertible and
$$
\operatorname{rank}(L-L')=s.
$$
\end{lemma}

\begin{proof}
Extend the rows of $R$ to a basis of $F^n$,
obtaining an invertible matrix
$$
N_0=\begin{pmatrix}R_0\\R\end{pmatrix}.
$$
Set $E=(I_s\ \ 0)$ and, for $t\in F$, define
$$
N(t)=\begin{pmatrix}R_0+tE\\R\end{pmatrix}.
$$
Let $D$ be the matrix formed by the first $s$
columns of $L-R_0$. Consider
$$
f(t)=\det N(t),
\qquad
g(t)=\det(D-tI_s).
$$
The polynomial $f$ is nonzero because
$f(0)=\det N_0\neq0$, while $g$ has degree $s$
and leading coefficient $(-1)^s$.
Since $F$ is infinite, there exists $t_0\in F$
such that
$$
f(t_0)g(t_0)\neq0.
$$
Taking $L'=R_0+t_0E$, we obtain an invertible
matrix $N=N(t_0)$. Moreover, the first $s$ columns
of $L-L'$ form the invertible matrix $D-t_0I_s$.
Thus $\operatorname{rank}(L-L')=s$.
\end{proof}
\section{Non-Surjectivity for large exponents}
\label{sec: not surjective}
In \cite{SAINI2027256}, the authors proved that, when $A_1$ is invertible, the inequality $n>k(r_0-1)$, where $r_0=\operatorname{nullity}(A_2)$, is necessary for the surjectivity of $\omega(x,y)=A_1x^k+A_2y^k$ and is also sufficient for $n=3,4$. We show below that this condition is not sufficient in general by constructing counterexamples in every order $n\geq6$.

\begin{proposition}\label{prop:nullity-two-counterexample}
Let $F$ be an algebraically closed field of characteristic
zero, let $n\geq6$, and let
$A\in M_{(n-2)\times n}(F)$ have rank $n-2$. Set
$$
B=\begin{pmatrix}A\\0_{2\times n}\end{pmatrix}.
$$
Then the map
$
\Phi(X,Y)=X^{n-1}+BY^{n-1}
$
is not surjective.
\end{proposition}

\begin{proof}
Set $p=n-2$.

We first observe that if $Z^{n-1}$ has a nontrivial
Jordan block at zero, then $Z$ is similar to $J_{0,n}$, see \cref{lem:Miller}.
Indeed, every zero Jordan block of size at most $n-1$
vanishes upon taking its $(n-1)$-st power. Consequently,
$$
\operatorname{rank}(Z^{n-1})=1.
$$

Since $\operatorname{rank}A=p\geq4$, choose
$G\in\mathrm{GL}_p(F)$ such that the first two rows
$v_1,v_2$ of $GA$ are linearly independent modulo
$Fe_n^\top$. Denote the remaining rows by $v_3,\dots,v_p$,
and define
$$
R=G^{-1}
\begin{pmatrix}
0\\
0\\
v_1\\
v_2\\
0_{(p-4)\times n}
\end{pmatrix},
\qquad
C=\begin{pmatrix}R\\e_n^\top\\0\end{pmatrix}.
$$
The block $0_{(p-4)\times n}$ in the definition of $R$ is
omitted when $p=4$. We show that $C$ is not in
the image of $\Phi$.

Suppose that
$$
C=M+BN,
\qquad M=X^{n-1},\quad N=Y^{n-1}.
$$
The last two rows give
$$
e_{n-1}^\top M=e_n^\top ,\qquad e_n^\top M=0.
$$
Thus $M$ has a nontrivial Jordan block at zero, and hence $\operatorname{rank}M=1$.
Since its $(n-1)$-st row is $e_n$, we may write
$$
M=\begin{pmatrix}a e_n^\top\\e_n^\top\\0\end{pmatrix}
$$
for some column vector $a\in F^p$.

Set $\alpha=Ga$. Multiplying the top-row equation
$R=a e_n^\top+AN$ on the left by $G$ gives
$$
\begin{aligned}
v_1N&=-\alpha_1e_n^\top,
&
v_2N&=-\alpha_2e_n^\top,\\
v_3N&=v_1-\alpha_3e_n^\top,
&
v_4N&=v_2-\alpha_4e_n^\top.
\end{aligned}
$$
The vectors $\{v_3N,v_4N\}$ are linearly independent because
$v_1,v_2$ are independent modulo $Fe_n^\top$. Hence
$\operatorname{rank}N\geq2$.

We next find a nonzero vector
$z\in\operatorname{span}(N)$ satisfying $zN=0$.
If $(\alpha_1,\alpha_2)\neq(0,0)$, the displayed
identities imply that $e_n^\top,v_1,v_2$ belong to
$\operatorname{span}(N)$. Therefore
$$
z=\alpha_2v_1-\alpha_1v_2
$$
is nonzero, belongs to $\operatorname{span}(N)$,
and satisfies $zN=0$.

If $\alpha_1=\alpha_2=0$, choose
$(b_1,b_2)\neq(0,0)$ such that
$b_1\alpha_3+b_2\alpha_4=0$. Then
$$
z=b_1(v_1-\alpha_3e_n^\top)+b_2(v_2-\alpha_4e_n^\top)
  =b_1v_1+b_2v_2
$$
is a nonzero element of $\operatorname{span}(N)$,
and $zN=0$.

In either case, there exists a row vector $w$ such that
$wN=z\neq0$ and $wN^2=0$. Thus $N$ has a nontrivial
Jordan block at zero. Since $N=Y^{n-1}$, the initial
observation gives $\operatorname{rank}N=1$, a contradiction.
\end{proof}
As a corollary, we obtain the following:
\begin{corollary}
\label{coro: counterexample}
Let $n\geq6$, and let $B\in M_n(F)$ have nullity two.
Then the map
$$
(X,Y)\longmapsto X^{n-1}+BY^{n-1}
$$
is not surjective.
\end{corollary}
\begin{proof}
By similarity, $B$ may be replaced by
$\begin{pmatrix}A\\0_{2\times n}\end{pmatrix}$,
where $A$ has rank $n-2$. The result follows from
Proposition~\ref{prop:nullity-two-counterexample}.
\end{proof}

\textbf{Proof of \cref{thm:large-exponents}}
Let $r_0=\operatorname{nullity}(A_2)$. If $r_0\leq1$, surjectivity follows from \cite[Theorem~1.2(1)]{SAINI2027256}. Suppose that $r_0\geq2$.
If $k\geq n$, then $k(r_0-1)\geq n$, so $\omega$ is not surjective by \cite[Theorem~1.2(2)]{SAINI2027256}.

It remains to consider $k=n-1$. If $r_0\geq3$, then $k(r_0-1)\geq2(n-1)\geq n$, and the same necessary condition gives non-surjectivity.
If $r_0=2$, set $B=A_1^{-1}A_2$.
Since $\operatorname{nullity}(B)=2$,
Corollary~\ref{coro: counterexample} shows that
$(X,Y)\mapsto X^{n-1}+BY^{n-1}$ is not surjective.
Multiplication by $A_1$ gives the required conclusion.

\section{The case \texorpdfstring{$n=5$}{n=5}}
\label{sec:n5}

In this section, we prove that the necessary condition
obtained in \cite{SAINI2027256} is also sufficient
in order five.

\begin{proposition}
    Let $B\in M_5(F)$ be a non-zero element and $k\geq 2$. Define the map 
    \begin{align*}
        \omega:M_5(F)\times M_5(F)\longrightarrow M_5(F), \qquad \omega(x,y)= x^k+By^k.
    \end{align*} The map $\omega$ is surjective if and only if $5>k(r_0-1)$.
\end{proposition}
\begin{proof}
The nullity of $B$ is denoted by $r_0$. For $r_0\leq1$, surjectivity follows from
\cite[Theorem~1.2(1)]{SAINI2027256}. For $r_0\geq2$, the inequality $5>k(r_0-1)$ is necessary by \cite[Theorem~1.2(2)]{SAINI2027256}. Thus, it remains
to prove that $\omega$ is surjective for $2\leq r_0\leq 3$ and $2\leq k\leq 4$ such that $k(r_0-1)<5$, which is precisely,
$$
(r_0,k)\in\{(2,2),(3,2),(2,3),(2,4)\}.
$$

Denote $s=5-r_0$. By the reduction in \cref{sec: Prelims},
we may assume that
$$
B=\begin{pmatrix}\widehat J\\0\end{pmatrix},
\qquad
\widehat J\in M_{s\times5}(F),
\qquad
\operatorname{rank}\widehat J=s.
$$
Let $C\in M_5(F)$ be arbitrary element such that
$$
C=\begin{pmatrix}L_1\\L_2\end{pmatrix},
\qquad
L_1\in M_{s\times5}(F) \text{ and } L_2\in M_{r_0\times5}(F).
$$

\textbf{Case 1: $k=2$.} Since $r_0\leq3$, every possible rank of $L_2$ satisfies one of the hypotheses of Lemma~\ref{lem: square completion}. Hence there exists $L_1'\in M_{s\times 5}(F)$ such that
$$
\begin{pmatrix}L_1'\\L_2\end{pmatrix}=X^2
$$
for some $X\in M_5(F)$. Set $H=L_1-L_1'$.
Choose $U\in\mathrm{GL}_5(F)$ satisfying
$$
\widehat J U=(0_{s\times r_0}\ \ I_s)=E.
$$
Since $s\leq3$, again by Lemma~\ref{lem: square completion}, there exists
$Z\in M_5(F)$ such that $EZ^2=HU$.
Taking $Y=UZU^{-1}$, we obtain
$$
\widehat JY^2
=\widehat JUZ^2U^{-1}
=EZ^2U^{-1}
=H.
$$
Therefore $X^2+BY^2=C$.

\textbf{Case 2: $r_0=2$, $k\in\{3,4\}$.} By Lemma~\ref{lem:residual lifting}, it suffices to
construct a $k$-th power matrix 
$N=\begin{pmatrix}L_1'\\L_2\end{pmatrix}$ such that the rank of 
$L_1-L_1'$ is at most one or $s$.
Indeed, writing $N=X^k$, the lemma gives
$Y\in M_n(F)$ satisfying
$\widehat JY^k=L_1-L_1'$, and therefore
$X^k+BY^k=C$.

 Set $E=(I_3\ \ 0)$. If $\operatorname{rank}(L_2)=2$, then
Lemma~\ref{lem:invertible-completion} gives
$L_1'\in M_{s\times n}(F)$ such that
$$
N=\begin{pmatrix}L_1'\\L_2\end{pmatrix}
$$
is invertible and
$\operatorname{rank}(L_1-L_1')=s$.
Thus $N$ is a $k$-th power with the required rank. 

If $L_2=0$, take $N(t)=tI_3\oplus0_2$. Clearly, $N=X^k$ for $X=\sqrt[k]{t}I_3\oplus0_2$. The determinant of the
first three columns of $L_1-tE$ has leading term
$-t^3$, so a suitable nonzero $t_0$ makes the
matrix $C-N$ have rank three.

It remains to consider $\operatorname{rank}L_2=1$. Recall that 
$
L_2=\begin{pmatrix}C_{21}^P&C_{22}^P\end{pmatrix}.
$
If $\operatorname{rank}C_{21}^P=1$, then $C\in \operatorname{Image}(\omega)$ by \cite[Proposition~4.5]{SAINI2027256} as $\dim (W_2)= \operatorname{rank}(C_{21}^P)$. We may therefore assume that $C_{21}^P=0$ and $\operatorname{rank}C_{22}^P=1$.

Let $Q\in GL_2(F)$ such that $QC_{22}^PQ^{-1}$ is the Jordan canonical form of $C_{22}^P$. Simultaneous conjugation by $S:=I_3\oplus Q$ preserves the zero lower rows
of $B$ and the rank of $\widehat J$. Putting
$C_{22}^P$ in Jordan canonical form, reduces $L_2$ to one of
$$
\begin{pmatrix}\mu e_4^\top\\0\end{pmatrix},
\qquad
\begin{pmatrix}e_5^\top\\0\end{pmatrix},
\qquad \mu\neq0.
$$

For the first form, take
$
N(t)=tI_3\oplus(\mu)\oplus0
$, which is a $k$-th power. A suitable nonzero $t_0$ makes $H=L_1-t_0E$ have rank three.

For the second form, if $\operatorname{rank}L_1\neq2$,
take $N=E_{45}$. By Lemma~\ref{lem: rank-one}, $N$ is a $k$-th power, and $C^S-N$ is of the form $\begin{pmatrix}
    L_1\\0_{2\times 5}
\end{pmatrix},$ having rank at most one or three.

Finally, suppose $\operatorname{rank}L_1=2$. Let
$c_1,c_2,c_3$ be its rows, and set
$$
W_1=\operatorname{span}\{c_1,c_2,c_3\}.
$$
If $e_5^\top\in W_1$, choose a basis
$\{e_5^\top,w\}$ of $W_1$ and write
$
c_i=a_i e_5^\top+b_iw $ for $1\leq i\leq3$.
Take
$
N=
\begin{pmatrix}a_1\\a_2\\a_3\\1\\0\end{pmatrix}
e_5^\top$. Then $N$ and $C-N$ both have rank at most one.

If $e_5^\top\notin W_1$, choose independent
rows $c_j,c_\ell$, and let $i$ be the remaining index.
Take
$$
N=(e_4+e_i)e_5^\top.
$$
The first three rows of $C-N$ are $c_j,c_\ell$ and
$c_i-e_5^\top$, in their original order.
They are linearly independent because
$e_5^\top\notin W_1$. Hence the matrix $C-N$
has rank three.

In both subcases, $N$ has rank one and is a $k$-th
power by Lemma~\ref{lem: rank-one}. Undoing the simultaneous conjugations completes the proof.
\end{proof}
As a corollary, we present the proof of Theorem~\ref{thm:n5-n6}(1):
\subsection{Proof of \texorpdfstring{\cref{thm:n5-n6}(1)}{Theorem (1)}}
\label{proof: thm n5} Since $A_1$ is invertible, let $B=A_1^{-1}A_2$. The $\operatorname{nullity}(B)=\operatorname{nullity}(A_2)$. Therefore, if $C=x^k+By^k$, then $A_1C=A_1x^k+A_2y^k$. Thus, $\omega_B$ is surjective if and only if $\omega$ is surjective. Hence, $\omega$ is surjective if and only if $5>k(r_0-1)$.

\section{The case \texorpdfstring{$n=6$}{n=6}}
\label{sec:n6}

As proved in Section~\ref{sec: not surjective}, the condition $n>k(r-1)$ is necessary but not sufficient for surjectivity. Indeed, for $n=6$, it holds when
$r=2$ and $k=5$, whereas the corresponding map is not
surjective, as shown in Proposition~\ref{prop:nullity-two-counterexample}.
In this section, we prove that for $n=6$ this is the only exception and obtain
the following classification.

We first show that every square-zero nilpotent matrix of rank two in order six has a $ k$-th root for $2\leq k\leq 4$.
\begin{lemma}\label{lem:n6-rank-two-nilpotent}
Let $N\in M_6(F)$ such that $\operatorname{rank}(N)=2$ and $N^2=0$. Then there exist $X\in M_6(F)$ such that $N=X^k$ for $k=2,3,4$.
\end{lemma}
\begin{proof}
The Jordan canonical form of $N$ is $J_{0,2}\oplus J_{0,2}\oplus0_2$. For $k=2$, take $X=(J_{0,4}\oplus0_2)$, for $k=3$ take $X=(J_{0,5}\oplus0)$ and $X=J_{0,6}$ for $k=4$. In each case, by \cref{lem:Miller}, $X^k$ is similar to $N$. 
\end{proof}
\begin{proposition}
Let $B\in M_6(F)$ be a non-zero element, and $k\geq 2$ be an integer. The map $\omega(x,y)=x^k+By^k$ is surjective on $M_6(F)$ if and only if $6>k(r_0-1)$ and $(r_0,k)\neq (2,5)$, where $r_0=\operatorname{nullity}(B)$.
\end{proposition}
\begin{proof}
For $r_0\leq1$, surjectivity follows from
\cite[Theorem~1.2(1)]{SAINI2027256}. For $r_0\geq2$,
the condition $6>k(r_0-1)$ is necessary by
\cite[Theorem~1.2(2)]{SAINI2027256}. The pairs satisfying
this condition are precisely
$$
(r_0,k)\in\{(2,2),(2,3),(2,4),(2,5),(3,2)\}.
$$
For $(r_0,k)=(2,5)$, the map $\omega$ is not surjective by
Proposition~\ref{prop:nullity-two-counterexample}. Thus, we prove the surjectivity for the remaining pairs.

Let $s=6-r_0$. By the reduction in \cref{sec: Prelims}, we may assume that $B=\begin{pmatrix}\widehat J\\0\end{pmatrix}, \text{ where }
\widehat J\in M_{s\times6}(F)\text{ and }
\operatorname{rank}\widehat J=s$.

Let $C\in M_6(F)$ be arbitrary, and we write $C=\begin{pmatrix}L_1\\L_2\end{pmatrix}$, where $L_1\in M_{s\times6}(F)$ and $
L_2\in M_{r_0\times6}(F)$.

\textbf{Case 1: $(r_0,k)=(3,2)$.} By Lemma~\ref{lem: square completion}, there exists
$X\in M_6(F)$ such that $X^2=\begin{pmatrix}L_1'\\L_2\end{pmatrix}$ for some $L_1'\in M_{3\times6}(F)$. Set $H=L_1-L_1'$.
Choose $U\in\mathrm{GL}_6(F)$ satisfying
$$
\widehat J U=(0_{3\times3}\ \ I_3)=E.
$$
Since $\operatorname{rank}(HU)\leq3$,
by Lemma~\ref{lem: square completion}, there exists $Z\in M_6(F)$ such that $EZ^2=HU$.
Taking $Y=UZU^{-1}$ gives
$$
\widehat JY^2=\widehat JUZ^2U^{-1}=H,
$$
and hence $X^2+BY^2=C$.

\textbf{Case 2: $r_0=2$, $2\leq k\leq4$.} By Lemma~\ref{lem:residual lifting}, it suffices to
construct a $k$-th power
$$
N=\begin{pmatrix}L_1'\\L_2\end{pmatrix}
$$
such that
$$
\operatorname{rank}(L_1-L_1')\in\{0,1,s\}.
$$
Indeed, writing $N=X^k$, the lemma gives
$Y\in M_n(F)$ satisfying
$\widehat JY^k=L_1-L_1'$, and therefore
$X^k+BY^k=C$.

 Set $E=(I_4\ \ 0)$. If $\operatorname{rank}(L_2)=2$, then
Lemma~\ref{lem:invertible-completion} gives
$L_1'\in M_{s\times n}(F)$ such that
$$
N=\begin{pmatrix}L_1'\\L_2\end{pmatrix}
$$
is invertible and
$\operatorname{rank}(L_1-L_1')=s$.
Thus $N$ is a $k$-th power with the required rank.

If $L_2=0$, take $N(t)=(tI_4\oplus0_2)= (\sqrt[k]{t}I_4\oplus0_2)^k$. A suitable nonzero $t_0$ makes $L_1-t_0E$ have rank four.

It remains to consider $\operatorname{rank}(L_2)=1$.
Recall \cref{sec: Prelims}, $L_2=\begin{pmatrix}C_{21}^P&C_{22}^P\end{pmatrix}$.

If $\operatorname{rank}C_{21}^P=1$, then $C$ belongs to the image by \cite[Proposition~4.5]{SAINI2027256}.
We may therefore assume that $C_{21}^P=0$ and
$\operatorname{rank}C_{22}^P=1$.

A simultaneous conjugation by $I_4\oplus Q$, with
$Q\in\mathrm{GL}_2(F)$, preserves the zero lower rows
of $B$ and the rank of $\widehat J$. Therefore, consider
$C_{22}^P$ in Jordan canonical form and thus $L_2$ reduces to one of
$$
\begin{pmatrix}\mu e_5^\top\\0\end{pmatrix},
\qquad
\begin{pmatrix}e_6^\top\\0\end{pmatrix},
\qquad \mu\neq0.
$$

For the first form, take $N(t)=tI_4\oplus(\mu)\oplus0=X^k$ for $X= \sqrt[k]{t}I_4\oplus(\sqrt[k]{\mu})\oplus0$. We can choose $t_0$ such that $H=L_1-t_0E$ has rank four.

Now, we consider the second form, $C=\begin{pmatrix}L_1\\e_6^\top\\0\end{pmatrix}$.
If $\operatorname{rank}L_1\in\{0,1,4\}$, take $N=E_{56}$, which is a $k$-th power by Lemma~\ref{lem: rank-one} and $C-N=\begin{pmatrix}
    L_1\\0
\end{pmatrix}$.

Suppose that $\operatorname{rank}L_1=2$. If $e_6^\top\in W_1$, choose a basis
$\{e_6^\top,w\}$ of $W_1$ and write the rows of $L_1$ as $c_i=a_i e_6^\top+b_iw$ for $1\leq i\leq 4$. Then 
$$\begin{pmatrix}
    L_1\\L_2
\end{pmatrix}= \begin{pmatrix}a_1\\a_2\\a_3\\a_4\\1\\0\end{pmatrix}
e_6^\top + \begin{pmatrix}b_1\\b_2\\b_3\\b_4\\0\\0\end{pmatrix}w=N+M.$$ Thus, $C$ is the sum of two rank-one matrices. The matrix $N=X^k$ for some $X$ by \cref{lem: rank-one} and $C-N$ is again a rank-one matrix.

If $e_6^\top\notin W_1$. Since $\operatorname{rank}(L_1)=2$, there exist $Q\in GL_4(F)$ such that $QL_1=\begin{pmatrix}
    u_1\\u_2\\0\\0
\end{pmatrix}$. Let $S=Q\oplus I_2$. Simultaneous conjugation of $C=x^k+Jy^k$
by $S$ gives $$\begin{pmatrix}
    QL_1\\ e_6^\top\\ 0
\end{pmatrix}=x^k+ \begin{pmatrix}
    Q\widehat J S^{-1}\\0_{2\times 6}
\end{pmatrix}y^k$$ where $QL_1S^{-1}=\begin{pmatrix}
    u_1S^{-1}\\u_2S^{-1}\\0\\0
\end{pmatrix}$.
The conjugation by $S$ preserves the last two rows of $C$, the form of $B$,
and the condition $e_6^\top\notin W_1$.

Consider $V=\operatorname{span}\{e_1^\top,e_2^\top, e_3^\top-e_5^\top\}$. Here $W_1$ denotes the row space of the transformed top block $QL_1S^{-1}$. Both $V$ and $W_1+Fe_6^\top$ have dimension three. However, $e_6^\top\in W_1+Fe_6^\top$ whereas $e_6^\top\notin V$. Hence these subspaces are distinct, and we may choose $w\in V\setminus(W_1+Fe_6^\top)$. Let $N=\begin{pmatrix}
0, 0, e_6^\top, w, e_6^\top, 0
\end{pmatrix}^\top$. Since $w\in V$, we have $w=(a,b,c,0,-c,0)$, we obtain $wN=0$. Also $e_6^\top N=0$, so $N^2=0$.
The row space of $N$ is spanned by $ w$ and $ e_6^\top$, so $\operatorname{rank}N=2$. By Lemma~\ref{lem:n6-rank-two-nilpotent}, $N$ is a
$k$-th power. 

Thus, we have
$$
C^S-N=
\begin{pmatrix}
u_1S^{-1}\\
u_2S^{-1}\\
-e_6^\top\\
-w\\
0_{2\times6}
\end{pmatrix}.
$$
The first two rows span $W_1$. Since
$e_6^\top\notin W_1$ and
$w\notin W_1+Fe_6^\top$, the first four rows are
linearly independent. Thus
$\operatorname{rank}(C^S-N)=4$.

Finally, suppose that $\operatorname{rank}L_1=3$.
If $e_6^\top\notin W_1$, choose three independent
rows of $L_1$, and let $i$ be the remaining index.
Taking $N=(e_5+e_i)e_6^\top$ gives a rank-one $k$-th power such that $C-N$ has rank four.

If $e_6^\top\in W_1$, a similar construction as above allows us to assume
that $L_1=\begin{pmatrix}u_1\\u_2\\u_3\\0\end{pmatrix}$ where $u_1,u_2,u_3$ are independent.
Let $V=\operatorname{span}\{e_1^\top,e_2^\top,e_3^\top\}$. Since $V$ and $W_1$ have dimension three and $e_6^\top\in W_1\setminus V$, choose $w\in V\setminus W_1$.
Define $N=\begin{pmatrix}0,0,0,w,e_6^\top,0\end{pmatrix}^\top$. Now, $wN=0$ and $e_6^\top N=0$, and thus $N^2=0$. Moreover, $\operatorname{rank}N=2$, so $N$ is a
$k$-th power by Lemma~\ref{lem:n6-rank-two-nilpotent}. The matrix $C-N$ has rank four, with the last two rows being zero.

In every case, we have constructed $N=X^k$
with the prescribed last two rows and $\operatorname{rank}(L_1-L_1')\in\{0,1,4\}$.

Lemma~\ref{lem:residual lifting} therefore gives
$Y\in M_6(F)$ satisfying
$\widehat JY^k=L_1-L_1'$.
Hence $X^k+BY^k=C$.
\end{proof}
As a corollary, we obtain:
\subsection{Proof of \texorpdfstring{\cref{thm:n5-n6}(2)}{Theorem (2)}}
\label{proof: n6}
 Since $A_1$ is invertible, let $B=A_1^{-1}A_2$. The $\operatorname{nullity}(B)=\operatorname{nullity}(A_2)$. Therefore, if $C=x^k+By^k$, then $A_1C=A_1x^k+A_2y^k$. Hence, $\omega$ is surjective if and only if
$6>k(r_0-1)$ and
$(r_0,k)\neq(2,5)$.

\section{Concluding remarks}
\label{sec: concluding}
The results above complete the surjectivity classification
in orders five and six and show that the necessary
condition $n>k(r_0-1)$ is not sufficient in any order
$n\geq6$. For coefficients of nullity two, the case
$k=n-1$ is settled by the general obstruction.
This leaves the following question for smaller exponents.

\begin{question}
Let $n\geq7$ and $2\leq k\leq n-2$.
For which matrices $B\in M_n(F)$ of nullity two
is the map
\[
\omega_B\colon M_n(F)\times M_n(F)\longrightarrow M_n(F),
\qquad
\omega_B(X,Y)=X^k+BY^k,
\]
surjective?
\end{question}
\printbibliography

@article{SAINI2027256,
  author  = {Saini, Prachi and Singh, Anupam},
  title   = {Polynomial maps with constants on matrix algebra},
  journal = {Linear Algebra and its Applications},
  volume  = {752},
  pages   = {256--277},
  year    = {2027},
  issn    = {0024-3795},
  doi     = {10.1016/j.laa.2026.08.016},
  url     = {https://www.sciencedirect.com/science/article/pii/S0024379526003356}
}

@article {PolynomialswithconstantsM,
    AUTHOR = {Panja, Saikat and Saini, Prachi and Singh, Anupam},
     TITLE = {Images of polynomial maps with constants},
   JOURNAL = {Mathematika},
  FJOURNAL = {Mathematika. A Journal of Pure and Applied Mathematics},
    VOLUME = {71},
      YEAR = {2025},
   %Number = {3},
     PAGES = {Paper No. e70031 (17 pages)},
     %%ISSN = {0025-5793,2041-7942},
   MRCLASS = {16S50 (11P05)},
  MRNUMBER = {4922424},
MRREVIEWER = {Omar\ Kchit},
       DOI = {10.1112/mtk.70031},
       URL = {https://doi.org/10.1112/mtk.70031},
}

@article {PolynomialsonSplitOctonion,
    AUTHOR = {Panja, Saikat and Saini, Prachi and Singh, Anupam},
     TITLE = {Polynomial maps with constants on split octonion algebras},
   JOURNAL = {Comm. Algebra},
  FJOURNAL = {Communications in Algebra},
    VOLUME = {54},
      YEAR = {2026},
   %Number = {2},
     PAGES = {579--594},
     %%ISSN = {0092-7872,1532-4125},
   MRCLASS = {17A75 (11P05 14M99 17A36 20B27)},
  MRNUMBER = {4996673},
       DOI = {10.1080/00927872.2025.2531559},
       URL = {https://doi.org/10.1080/00927872.2025.2531559},
}

@article {WordmapsandWaringproblemLArsen,
    AUTHOR = {Larsen, Michael and Shalev, Aner},
     TITLE = {Word maps and {W}aring type problems},
   JOURNAL = {J. Amer. Math. Soc.},
  FJOURNAL = {Journal of the American Mathematical Society},
    VOLUME = {22},
      YEAR = {2009},
   % NUMBER = {2},
     PAGES = {437--466},
%      %ISSN = {0894-0347,1088-6834},
   MRCLASS = {20D06 (14G15 20G40 20P05)},
  MRNUMBER = {2476780},
MRREVIEWER = {Martin\ W.\ Liebeck},
       DOI = {10.1090/S0894-0347-08-00615-2},
       URL = {https://doi.org/10.1090/S0894-0347-08-00615-2},
}

@article {Wordmapswithconstants,
    AUTHOR = {Gordeev, Nikolai and Kunyavski\u i, Boris and Plotkin, Eugene},
     TITLE = {Word maps, word maps with constants and representation
              varieties of one-relator groups},
   JOURNAL = {J. Algebra},
  FJOURNAL = {Journal of Algebra},
    VOLUME = {500},
      YEAR = {2018},
     PAGES = {390--424},
      %ISSN = {0021-8693,1090-266X},
   MRCLASS = {20G15 (20F10 20F70)},
  MRNUMBER = {3765462},
MRREVIEWER = {Alla\ S.\ Detinko},
       DOI = {10.1016/j.jalgebra.2017.03.016},
       URL = {https://doi.org/10.1016/j.jalgebra.2017.03.016},
}

@article {GnutovWordmapswithcosntants,
    AUTHOR = {Gnutov, F. A.},
     TITLE = {On the image of a word map with constants of a simple
              algebraic group {II}},
   JOURNAL = {Zap. Nauchn. Sem. S.-Peterburg. Otdel. Mat. Inst. Steklov.
              (POMI)},
  FJOURNAL = {Rossi\u iskaya Akademiya Nauk. Sankt-Peterburgskoe Otdelenie.
              Matematicheski\u i\ Institut im. V. A. Steklova. Zapiski
              Nauchnykh Seminarov (POMI)},
    VOLUME = {492},
      YEAR = {2020},
     PAGES = {75--93},
     %% %ISSN = {0373-2703},
   MRCLASS = {20G15 (20F10)},
  MRNUMBER = {4223709},
MRREVIEWER = {Denis\ E.\ Serbin},
Note={MR4223709}
}

@misc{bresar2026waring,
  title         = {Waring Problems Across Algebra},
  author        = {Matej Bre{\v{s}}ar and Consuelo Mart{\'i}nez},
  year          = {2026},
  eprint        = {2603.07033},
  archivePrefix = {arXiv},
 % primaryClass  = {math.RA},
  url           = {https://arxiv.org/abs/2603.07033}
}

@article {BresarSemrlwaringproblem,
    AUTHOR = {Bre\v sar, Matej and \v Semrl, Peter},
     TITLE = {The {W}aring problem for matrix algebras, {II}},
   JOURNAL = {Bull. Lond. Math. Soc.},
  FJOURNAL = {Bulletin of the London Mathematical Society},
    VOLUME = {55},
      YEAR = {2023},
    NUMBER = {4},
     PAGES = {1880--1889},
      %ISSN = {0024-6093,1469-2120},
   MRCLASS = {16R10 (16S50)},
  MRNUMBER = {4623691},
MRREVIEWER = {D.\ Diniz},
       DOI = {10.1112/blms.12825},
       URL = {https://doi.org/10.1112/blms.12825},
}

@article {MR4592319,
    AUTHOR = {Panja, Saikat and Prasad, Sachchidanand},
     TITLE = {The image of polynomials and {W}aring type problems on upper
              triangular matrix algebras},
   JOURNAL = {J. Algebra},
  FJOURNAL = {Journal of Algebra},
    VOLUME = {631},
      YEAR = {2023},
     PAGES = {148--193},
      %ISSN = {0021-8693,1090-266X},
   MRCLASS = {16S50 (16R10)},
  MRNUMBER = {4592319},
MRREVIEWER = {Pedro\ Souza\ Fagundes},
       DOI = {10.1016/j.jalgebra.2023.04.027},
       URL = {https://doi.org/10.1016/j.jalgebra.2023.04.027},
}

@article {MR4526834,
    AUTHOR = {Gargate, Ivan Gonzales and de Mello, Thiago Castilho},
     TITLE = {Images of multilinear polynomials on {$n$} {$\times$} {$n$}
              upper triangular matrices over infinite fields},
   JOURNAL = {Israel J. Math.},
  FJOURNAL = {Israel Journal of Mathematics},
    VOLUME = {252},
      YEAR = {2022},
    NUMBER = {1},
     PAGES = {337--354},
      %ISSN = {0021-2172,1565-8511},
   MRCLASS = {16S50 (15A69 16R10)},
  MRNUMBER = {4526834},
MRREVIEWER = {Ronald\ Ismael\ Quispe Urure},
       DOI = {10.1007/s11856-022-2350-2},
       URL = {https://doi.org/10.1007/s11856-022-2350-2},
}

@misc{demello2026generalizedcommutatoroctonionalgebra,
      title={Generalized Commutator on Octonion Algebra}, 
      author={Thiago Castilho De Mello and Prachi Saini},
      year={2026},
      eprint={2610.04466},
      archivePrefix={arXiv},
      primaryClass={math.RA},
      url={https://arxiv.org/abs/2610.04466}, 
}

@article {MR4575412,
    AUTHOR = {Bre\v sar, Matej and \v Semrl, Peter},
     TITLE = {The {W}aring problem for matrix algebras},
   JOURNAL = {Israel J. Math.},
  FJOURNAL = {Israel Journal of Mathematics},
    VOLUME = {253},
      YEAR = {2023},
    NUMBER = {1},
     PAGES = {381--405},
      ISSN = {0021-2172,1565-8511},
   MRCLASS = {16R10 (16S50)},
  MRNUMBER = {4575412},
MRREVIEWER = {D.\ Diniz},
       DOI = {10.1007/s11856-022-2366-7},
       URL = {https://doi.org/10.1007/s11856-022-2366-7},
}
\end{document}